\documentclass[12pt,thmsa]{article}
\usepackage{amsmath}
\usepackage{graphicx}
\usepackage{amsfonts}
\usepackage{amssymb}
\usepackage{mathrsfs}

\newtheorem{prop}{Proposition}
\newtheorem{lem}{Lemma}
\newtheorem{thm}{Theorem}

\begin{document}

\begin{center}
{\Large \textbf{Kernel-based method for joint independence of functional variables}}

\bigskip

Terence Kevin  MANFOUMBI DJONGUET and  Guy Martial  NKIET 

\bigskip

\textsuperscript{}URMI, Universit\'{e} des Sciences et Techniques de Masuku,  Franceville, Gabon.

\bigskip

E-mail adresses : tkmpro95@gmail.com,    guymartial.nkiet@univ-masuku.org.

\bigskip
\end{center}

\noindent\textbf{Abstract.}  This work investigates the problem of testing whether $d$ functional random variables are jointly independent using a modified estimator of the $d$-variable Hilbert Schmidt Indepedence Criterion ($d$HSIC) which  generalizes  HSIC for the case  where $d \geq 2$. We then get asymptotic normality of this estimator both under joint independence hypothesis and under the alternative hypothesis.  A   simulation study shows good performance  of  the proposed test on finite sample.
\bigskip


\vskip 4mm

\noindent 1.   INTRODUCTION

\vskip 3mm

\noindent The question of relationships between variables has always been a major concern in statistical analysis. This is why methods providing answers to this type of problems, such as regression methods or independence tests, occupy  an important place in the statistical literature. However, there exists just a few works that deal with independence testing in the context of functional data analysis, despite the importance of this field over the last two decades.   Kokoszka et al. (2008) and  Aghoukeng Jiofack and Nkiet (2010) introduced  non-correlation tests  between two functional variables, whereas   G\'orecki et al. (2020),  Lai et al. (2021) and  Lauman et al. (2021)  recently  proposed tests for independence. Kernel-based methods, that use the distance between the kernel embeddings of probability measures in a reproducing kernel Hilbert space (RKHS), were also  introduced in Gretton et al. (2005, 2007). For the case where  more than two variables are considered,  pairwise independence tests can be obtained by perfoming the aforementioned tests  on each pair,  but  it is not sufficient for  some statistical problems such as those arising in   causal inference where models often assume the existence of jointly independent noise variables, or in independent component analysis. So, it can be more relevant to consider testing for joint independence of several random variables rather  than testing pairwise independence. Recently, Pfister et al. (2018) proposed    kernel-based methods  for testing for joint independence of $d$ random variables valued into metric spaces, with $d\geq 2$. For doing that, they introduced the $d$-variable Hilbert-Schmidt independence criterion ($d$HSIC) as a measure of joint independence that extends the Hilbert-Schmidt independence criterion (HSIC)  of Gretton et al. (2005) and contains it as a special case. For testing for joint independence  they considered, as test statistic, a consistent estimator of $d$HSIC and derived  its  asymptotic distribution under null hypothesis. However, as it is the case for most measures introduced in kernel-based methods such as maximal mean discrepancy (MMD) (Gretton et al. (2012)), generalized  maximal mean discrepancy (GMMD) (Balogoun et al. (2021,2022)) and HSIC, this asymptotic distribution     is an infinite sum of distributions and, therefore,  can not be used for performing the test.   That is why Pfister et al. (2018)  proposed three approaches   based on $d$HSIC: a permutation test, a bootstrap test and a gamma approximation of the limiting distribution. An alternative approach can be tackled as it was done   in  Makigusa and Naito (2020) and  in Balogoun et al. (2021) for the MMD and the GMMD. It consists  in constructing a test statistic  from   an  appropriate modification of  a naive  estimator   in order to yield asymptotic normality both under the  null hypothesis and under the alternative. In this paper,  such an approach is adopted for the $d$HSIC. It   allows us to propose a new   test for joint independence of  several  random variables  valued into  metric spaces, including functional variables. The rest of the paper is organized as follows.     The  $d$HSIC  is recalled in Section 2, and Section 3 is devoted to its estimation by a modification of a  naive estimator,  and to the main results.  A simulation study on functional data that allows to compare   the proposed test   to that of  Pfister  et al. (2018) is  given in Section 4. All the proofs are postponed in Section 5.
\vskip 3mm

\noindent 2.   $d$HSIC AND JOINT INDEPENDENCE PROPERTY

\vskip 3mm

\noindent For an integer $d\geq 2$, let  $X^{(1)},\cdots,X^{(d)}$ be   random variables defined on a probability space $(\Omega,\mathscr{A},P)$ and with values in compact metric spaces $\mathcal{X}_1,\cdots,\mathcal{X}_d$  respectively. Then, we consider the random variable $$X=(X^{(1)},\cdots,X^{(d)})$$ with values in $\mathcal{X}:=\mathcal{X}_1\times\mathcal{X}_2\times\cdots\times\mathcal{X}_d$. For $\ell\in\{1,\cdots,d\}$,  let us introduce  a  reproducing kernel Hilbert space (RKHS)    $\mathcal{H}_\ell$    of functions  from  $\mathcal{X}_\ell$  to $\mathbb{R}$  with associated   kernel    $K_\ell\,:\,\mathcal{X}_\ell^2\rightarrow\mathbb{R}$. Each $K_\ell$ is a   symmetric function  such that, for any $f\in \mathcal{H}_\ell$ and any $t\in\mathcal{X}_\ell$, one has $K_\ell(t,\cdot)\in\mathcal{X}_\ell$ and  $f(t)=\langle K_\ell(t,\cdot),f\rangle_{\mathcal{X}_\ell}$  (see Berlinet and Thomas-Agnan (2004)), where  $\langle \cdot , \cdot \rangle_{\mathcal{X}_\ell}$  denotes the inner product of $\mathcal{X}_\ell$.    Denoting by $f_1\otimes\cdots\otimes f_d$ the tensor products of functions defined by:
\[
\left(f_1\otimes\cdots\otimes f_d\right)\left(x_1,\cdots,x_d\right)=\prod_{\ell=1}^df_\ell(x_\ell),
\]
we consider $\mathcal{H}:=\mathcal{H}_1\otimes\cdots\otimes\mathcal{H}_d$ the completed of the vector space spanned by the set $\{f_1\otimes\cdots\otimes f_d\,/\,f_\ell\in\mathcal{H}_\ell,\,\ell=1,\cdots,d\}$. It is known that $\mathcal{H}$ is a RKHS associated with the kernel $K$ satisfying:
\begin{equation}\label{produit}
K\left(x,y\right)=\prod_{\ell=1}^dK_\ell(x_\ell,y_\ell)
\end{equation}
for any $x=(x_1,\cdots,x_d)\in\mathcal{X}$ and  $y=(y_1,\cdots,y_d)\in\mathcal{X}$. Throughout this paper, we make the following assumption:

\bigskip

\noindent $(\mathscr{A}_1):$  $\Vert  K_\ell\Vert_{\infty}:= \sup\limits_{(t,s)\in \mathcal{X}_\ell^2} K_\ell(t,s)<+ \infty $,\,\,\,$\ell=1,\cdots,d$;

\bigskip
\noindent it implies that  $\Vert  K\Vert_{\infty}<+ \infty$ and  ensures the existence of the kernel mean embeddings $m_\ell=\mathbb{E}\left(K_\ell(X^{(\ell)},\cdot)\right)$, $\ell=1,\cdots,d$, and $m=\mathbb{E}\left(K(X,\cdot)\right)$. The random variables  $X^{(1)},\cdots,X^{(d)}$ are said to be  jointly independent   if
$$ \mathbb{P}_{X} = \mathbb{P}_{X^{(1)}} \otimes \mathbb{P}_{X^{(2)}} \otimes \cdots \otimes\mathbb{P}_{X^{(d)}},$$
where $\mathbb{P}_{X}$ (resp. $\mathbb{P}_{X^{(j)}}$) denotes the distribution of $X$ (resp. $X^{(j)}$). For measuring this property, Pfister et al. (2018) introduced the $d$-variable  Hilbert-Schimdt independence criterion    (dHSIC)   defined as
\begin{equation}\label{dhsic}
d\textrm{HSIC}(X) 
 =\bigg \Vert    \mathbb{E}\left(K(X,\cdot)\right)  - \bigotimes_{\ell=1}^{d} m_\ell \bigg\Vert _\mathcal{H} ^2
=\bigg \Vert    \mathbb{E}\left(\prod_{\ell=1}^dK_\ell(X^{(\ell)},\cdot)\right)  - \bigotimes_{\ell=1}^{d} m_\ell \bigg\Vert _\mathcal{H} ^2,
\end{equation}
where   $\bigotimes_{\ell=1}^{d} m_\ell = m_1 \otimes m_2 \otimes \dots \otimes m_d$. It is known  that if the kernels $K_\ell$ are characteristic ones, then  $d$HSIC fully characterizes joint  independence since the equality  $\mathbb{P}_{X} = \mathbb{P}_{X^{(1)}} \otimes \mathbb{P}_{X^{(2)}} \otimes \cdots \otimes\mathbb{P}_{X^{(d)}}$ is equivalent to $d\textrm{HSIC}(X)=0$  (see Pfister et al. (2018)).  So, in order to test for joint independence, that is testing for the hypothesis $\mathscr{H}_0\,:\,\mathbb{P}_{X} = \mathbb{P}_{X^{(1)}} \otimes \mathbb{P}_{X^{(2)}} \otimes \cdots\otimes \mathbb{P}_{X^{(d)}}$ against $\mathscr{H}_1\,:\,\mathbb{P}_{X}\neq \mathbb{P}_{X^{(1)}} \otimes \mathbb{P}_{X^{(2)}} \otimes \cdots \otimes\mathbb{P}_{X^{(d)}}$, a consistent  estimator of  $d\textrm{HSIC}(X)$ can be used as test statistic. Let $\{X_i\}_{1\leq i\leq n}$ be an i.i.d. sample of $X$, with $X_i=(X^{(1)}_i,\cdots,X^{(d)}_i)$,   Pfister et al. (2018) introduced the statistic  $\widehat{d\textrm{HSIC}}_{n}$ given by
\begin{eqnarray}\label{pfister}
\widehat{d\textrm{HSIC}}_{n}
& =  & \frac{1}{n^2} \sum_{i=1}^{n} \sum_{j=1} ^n\bigg( \prod_{\ell=1}^{d} K_\ell (X_i ^{(\ell)},X_j ^{(\ell)}) \bigg) + \frac{1}{n^{2d}}  \sum_{i_1=1} ^n\cdots \sum_{i_{2d}=1} ^n \prod_{\ell=1}^{d}  K_\ell (X_{i_{2\ell-1}} ^{(\ell)},X_{i_{2\ell}} ^{(\ell)})\nonumber\\
&  &- \frac{2}{n^{d+1}}\sum_{i_1=1} ^n\cdots \sum_{i_{d+1}=1} ^n \prod_{\ell=1}^{d}  K_\ell (X_{i_{\ell}} ^{(\ell)},X_{i_{\ell+1}} ^{(\ell)})
\end{eqnarray}
if $n\geq 2d$, and $\widehat{d\textrm{HSIC}}_{n}=0$ if $n\in\{1,\cdots,2d-1\}$. They proved that, under $\mathscr{H}_0$, the sequence  $n\,\, \widehat{d\textrm{HSIC}}_{n}$ converges in distribution, as $n\rightarrow +\infty$, to 
$\binom{2d}{d}\sum_{m=1}^{+\infty}\lambda_m\,Z_m^2$, where $(\lambda_m)_{m\geq 1}$ is an appropriate sequence of positive real numbers and   $(Z_m)_{m\geq 1}$ is a sequence of independent standard normal random variables.  This limiting distribution cannot be used for performing the test because it is an infinite sum of distributions. That is why Pfister et al. (2018) introduced three approach for this test based on the above estimator: a permutation test, a boostrap analogue and a procedure based on a gamma approximation. We propose in this paper another estimator for which we can obtain asymptotic normality under the null hypothesis. The main advantage is that it permits to avoid the use of  permutation or bootstrap methods for achieving the test, so leading to a  faster procedure. 

\vskip 3mm

\noindent 3. MODIFIED ESTIMATOR  AND ASYMPTOTIC NORMALITY

\vskip 3mm

\noindent A  natural approach for estimating  $d\textrm{HSIC}(X) $ consists in   replacing in \eqref{dhsic} each expectation by its empirical counterpart. This leads to the naive estimator
\begin{eqnarray}\label{estimator}
 \widetilde{D} _{n}
& =&\bigg \Vert    \frac{1}{n} \sum_{i=1}^{n} K(X_i,\cdot)   - \bigotimes_{\ell=1}^{d} \widehat{m}_\ell \bigg\Vert _\mathcal{H} ^2\nonumber \\
& =  & \frac{1}{n^2} \sum_{i=1}^{n} \sum_{j=1} ^n\bigg( \prod_{\ell=1}^{d} K_\ell (X_i ^{(\ell)},X_j ^{(\ell)}) \bigg) + \frac{1}{n^{2d}}  \prod_{\ell=1}^{d} \bigg(  \sum_{i=1}^{n} \sum_{j=1} ^n  K_\ell (X_i ^{(\ell)},X_j ^{(\ell)}) \bigg)\nonumber\\
&  &- \frac{2}{n^{d+1}} \sum_{i=1}^{n}  \prod_{\ell=1}^{d}\bigg( \sum_{j=1}^{n} K_\ell (X_i ^{(\ell)},X_j ^{(\ell)})\bigg),
\end{eqnarray}
where $\widehat{m}_\ell = \frac{1}{n} \sum_{i=1}^{n}  K_\ell (X_i ^{(\ell)},\cdot)$.  Instead, we adopt an approach introduced in Makigusa and Naito (2020), and also  used in Balogoun et al. (2021), consisting in introducing a modification  in \eqref{estimator} in order to get   a test statistic for which asymptotic normality can be obtained both under $\mathscr{H}_0$ and  under  $\mathscr{H}_1$.  For $\gamma \in  \; ]0,1]$, let $(w_{i,n}(\gamma))_{1 \leq i \leq n}$ be a sequence of positive numbers satisfying :\\

\bigskip
\noindent $(\mathscr{A}_2):$ There exists  a strictly positive real number $ \tau $  and an  integer $ n_0 $ such that for all $n> n_0$:
$$n\left\vert\frac{1}{n}\sum_{i=1}^{n}w_{i,n}(\gamma)-1\right\vert\leq \tau.$$

\noindent $(\mathscr{A}_3):$ There exists  $C>0$ such that $\max\limits_{1\leq i \leq n}w_{i,n}(\gamma)<C$ for all  $n\in\mathbb{N}^\ast$ and all $ \gamma \in]0, 1]$.

\noindent$(\mathscr{A}_4):$ For any   $ \gamma \in]0, 1]$,
$\lim\limits_{n\rightarrow +\infty}\frac{1}{n}\sum_{i=1}^{n}w^2_{i,n}(\gamma)=w^2(\gamma)> 1$.

 \bigskip
\noindent An example of  such sequence was given in Ahmad (1993)  as   $w_{i,n}(\gamma)=1+(-1)^i\,\gamma$.  For this example, one has  $C=2$, $w^2(\gamma)=1+\gamma^2$ and $\tau$ is any positive real number. Another example is   $w_{i,n}(\gamma)=1+\sin(i\pi\gamma)$ which corresponds to $\tau=1/\vert\sin(\pi\gamma/2)\vert$, $C=2$ and $w^2(\gamma)=3/2$.

\bigskip

\noindent We then propose to estimate $d\textrm{HSIC}(X)$ by a modification of \eqref{estimator} given by
\begin{eqnarray}\label{modifiedestimator}
\widehat{\textrm{\textbf{D}}}_{\gamma, n}
& =  & \frac{1}{n^2} \sum_{i=1}^{n} \sum_{j=1} ^n\bigg( \prod_{\ell=1}^{d} K_\ell (X_i ^{(\ell)},X_j ^{(\ell)}) \bigg) + \frac{1}{n^{2d}}  \prod_{\ell=1}^{d} \bigg(  \sum_{i=1}^{n} \sum_{j=1} ^n  K_\ell (X_i ^{(\ell)},X_j ^{(\ell)}) \bigg)\nonumber\\
&  &- \frac{2}{n^{d+1}} \sum_{i=1}^{n}  w_{i,n}(\gamma) \prod_{\ell=1}^{d}\bigg( \sum_{j=1}^{n} K_\ell (X_i ^{(\ell)},X_j ^{(\ell)})\bigg),
\end{eqnarray}
and we take this estimator as test statistic for testing for $\mathscr{H}_0$ against $\mathscr{H}_1$. For achieving this test, the asymptotic distribution under $\mathscr{H}_0$ of the test statistic is needed. We will obtain an asymptotic normality result  which is usable both under  $\mathscr{H}_0$ and under  $\mathscr{H}_1$, so permitting to perform the test. For any $(k,q)\in\mathbb{N}^\ast$ such that $q\geq k$, we put   $\nu_{k}^q = \bigotimes_{\ell=k}^{q} m_\ell$. Then, considering the functions $\mathcal{U}$ and  $\mathcal{V}$  from $\mathcal{X}$ to $\mathbb{R}$ defined on any $x=\left(x^{(1)},\cdots,x^{(d)}\right)\in\mathcal{X}$ by 
\begin{eqnarray*}
\mathcal{U} (x) &=& \Big\langle  K_1(x^{(1)},\cdot)\otimes \nu_{2}^{d} \\
&& + \sum_{\ell=2}^{d-1}   \left(\nu_{1}^{\ell-1} \otimes K_\ell(x ^{(\ell)},\cdot)\otimes \nu_{\ell+1}^{d}\right) + \nu_{1}^{d-1}\otimes K_d(x ^{(d)},\cdot) -  d\,\nu_{1}^d, \nu_{1}^d- m  \Big \rangle_ \mathcal{H}\\
& & +\Big \langle	 K(x,\cdot) - m   , m   \Big \rangle_ \mathcal{H}  
\end{eqnarray*}
and
$
\mathcal{V} (x) =  \Big \langle	 K(x,\cdot) - m   , \nu_{1}^d\Big \rangle_ \mathcal{H}$, 
we have:

\bigskip

\begin{thm}{}\label{thm} Assume that ($\mathscr{A}_1$) to ($\mathscr{A}_4$) hold. Then as $n\rightarrow +\infty$, one has  $$ \sqrt{n}\Big (\widehat{\textit{\textbf{D}}}_{\gamma,n}-   d\textrm{HSIC}(X) \Big ) \stackrel{\mathscr D}{\longrightarrow} \mathcal N (0, \sigma ^2_\gamma),$$ where: 
\begin{equation*}
\begin{aligned}
\sigma ^2_\gamma   & =  4 Var\left(\mathcal{U}(X_1)\right)+4 w^2(\gamma)Var\left(\mathcal{V}(X_1)\right) -8Cov\left(\mathcal{U}(X_1),\mathcal{V}(X_1)\right).
\end{aligned}
\end{equation*}
\end{thm}
\bigskip
This theorem  gives asymptotic normality in the general case. The particular case where  $\mathscr{H}_0$ holds can then  be deduced. Indeed, in this case one has $\nu_1^d=\bigotimes_{\ell=1}^{d} m_\ell=m$ and, therefore, $\mathcal{U}=\mathcal{V}$. Since $d\textrm{HSIC}(X)=0$ it follows that $\sqrt{n}\widehat{\textit{\textbf{D}}}_{\gamma,n}  \stackrel{\mathscr D}{\longrightarrow} \mathcal N (0, \sigma ^2_\gamma)$, where: 
\begin{equation*}
\begin{aligned}
\sigma_\gamma ^2 & =    4 (w^2 (\gamma)- 1) Var\left(\mathcal{V} (X_1)\right).
\end{aligned}
\end{equation*}
This asymptotic variance in unknown since it depends on the distribution of $X$ which is unknown. So, in order to perform the test, we have to seek a consistent estimator of it. The following proposition gives such an estimator:

\bigskip
\begin{prop} \label{prp} Assume that ($\mathscr{A}_1$) to ($\mathscr{A}_4$) hold, and put
$\widehat{\sigma}_\gamma^2=4 (w(\gamma) ^2- 1)\widehat{\alpha}^2_n,$ where
\begin{eqnarray*}
\widehat{\alpha}_n^2=\frac{1}{n}\sum_{i=1}^{n}  \Bigg ( \frac{1}{n}\sum_{j=1}^{n}  \Bigg[\prod_{\ell=1}^{d}K_\ell(X^{(\ell)}_i,X^{(\ell)}_j)\Bigg]\Bigg)^2-\bigg(\frac{1}{n^2} \sum_{i=1}^{n} \sum_{j=1}^{n}  \Bigg[\prod_{\ell=1}^{d}K_\ell(X^{(\ell)}_i,X^{(\ell)}_j)\Bigg]\bigg)^2.
\end{eqnarray*}
Then, under $\mathscr{H}_0$, $\widehat{\sigma}_\gamma^2$ converges in probability to  $\sigma_\gamma ^2$   as $n\rightarrow +\infty$.
\end{prop}

\bigskip

\noindent We deduce from the previous  theorem and proposition that, under  the hypothesis  $\mathscr{H}_0$,   $\sqrt{n}\,\widehat{\sigma}_\gamma^{-1}\widehat{\textit{\textbf{D}}}_{\gamma,n}\stackrel{\mathscr D}{\longrightarrow} \mathcal N (0, 1)$ as $n\rightarrow +\infty$.  This allows to perform the test:   for a given significant level $\alpha \in ]0,1[$, one has to reject $\mathscr{H}_0$ if  $\big\vert \widehat{\textit{\textbf{D}}}_{\gamma,n}\big\vert>  n^{-1/2} \widehat{\sigma}_\gamma \Phi^{-1}(1-\alpha /2)$ where $\Phi$ is the cumulative distribution function of the standard normal distribution.

\bigskip

\noindent\textbf{Remark 1.}
This test can easily  be applied on functional data. Indeed, suppose that  each  $X_i^{(\ell)}$ is a  random functions belonging to $L^2([0,1])$ and observed on points $t_1^{(\ell)},\cdots,t_{r_\ell}^{(\ell)}$   of a   fine grid  in $[0,1]$ such that $t_1^{(\ell)}= 0$ and $t_{r_\ell}^{(\ell)}= 1$. 
If   gaussian kernels  are  used, one has
\[
K_\ell(X^{(\ell)}_i ,X^{(\ell)}_j )=\exp\left(-\eta_\ell^2\Vert X_i-X_j\Vert^2\right)=\exp\left(-\eta_\ell^2\int_0^1\left(X^{(\ell)}_i(t)-X^{(\ell)}_j(t)\right)^2\,dt\right),
\]
where $\eta_\ell>0$, and this term  can   be approximated  by using trapezoidal rule:
\begin{eqnarray}\label{approx}
K_\ell(X^{(\ell)}_i ,X^{(\ell)}_j )&\simeq&\exp\bigg(-\eta_\ell^2\sum_{m=1}^{r_\ell-1}\bigg\{\frac{t^{(\ell)}_{m+1}-t^{(\ell)}_m}{2}\bigg(\left(X^{(\ell)}_i(t^{(\ell)}_m)-X^{(\ell)}_j(t^{(\ell)}_m)\right)^2\nonumber\\
&&\hspace{3cm}+\left(X^{(\ell)}_i(t^{(\ell)}_{m+1})-X^{(\ell)}_j(t^{(\ell)}_{m+1})\right)^2\bigg)\bigg\}\bigg).
\end{eqnarray}
Then,   $ \widehat{\textit{\textbf{D}}}_{\gamma,n}$ and $\widehat{\alpha}^2_n$  are to be computed by using \eqref{approx}.

\vskip 3mm

\noindent 4. SIMULATIONS

\vskip 3mm

\noindent In this section, we present a  simulation study   made in order to evaluate   the finite sample performance of the proposed test and compare it to  the test  of Pfister  et al. (2018)  based on a permutation method using $\widehat{d\textrm{HSIC}}_{n}$.  For convenience, we denote our test as T1, and the test of  Pfister  et al. (2018)   as   T2.   We computed   empirical sizes and powers    through   Monte Carlo simulations.  We considered  the case where $d=3$,  $\mathcal{X}_1=\mathcal{X}_2=\mathcal{X}_3=\textrm{L}^2([0,1])$, and used the following  models for generating functional data:
\[
X^{(1)}(t)=\sqrt{2}\sum_{k=1}^{50}\alpha_k\,\cos(k\pi t),\,X^{(2)}(t)=\sqrt{2}\sum_{k=1}^{50}\beta_k\,\cos(k\pi t)+\lambda f\left(X^{(1)}(t)\right),
\]
\[
X^{(3)}(t)=\sqrt{2}\sum_{k=1}^{50}\xi_k\,\cos(k\pi t)+\lambda f\left(X^{(1)}(t)\right)+\lambda f\left(X^{(2)}(t)\right),
\]
 where the $\alpha_k$s, the $\beta_k$s and the $\xi_k$s are independent and i.i.d. having the uniform distribution on $[0,1]$, $f$ is a given function and $\lambda$ is a real number giving the level  of dependence between the functional variables. Indeed,    $\mathscr{H}_0$ holds in case $\lambda=0$, and the dependence level increases with $\lambda$. Empirical sizes and powers were computed on the basis of $200$ independent  replicates. For each of them, we generated   samples  of size   $n=100$ of the above  processes in discretized versions  on equispaced values $t_1,\cdots,t_{51}$  in  $[0,1]$, where $t_j=(j-1)/50$, $j=1,\cdots,51$.  For performing  our method, we took $\gamma=0.32$ and used the gaussian  kernels 
\[
K_1(x,y)=K_2(x,y)=K_3(x,y)=\exp\left(-150\int_0^1\left( x(t)-y(t)\right)^2dt\right).
\] 
The test statistics given in \eqref{pfister} and \eqref{modifiedestimator} were computed by   approximating   integrals involved in the above    kernels  by using the trapezoidal rule, as indicated in \eqref{approx}. The significance level was taken as $\alpha=0.05$. The method T2 was used with $100$ permutations. Table 1 reports the obtained results. They are    close to the nominal size for both  methods when  $\lambda=0$. For $\lambda=0.25$, T2 gives  slightly larger powers than T1, but   the differences  are    low. For $m=0.5,1$, the two methods give maximal  power, except for the case where $f(x)=\sin(x)$ for which the powers are surprisingly low.  This highlights the interest of the proposed test:  it is powerful enough and is fast compared to the method of Pfister et al. (2018)   based on  permutations  wich leads to   very high computation times.

\begin{table}
{\setlength{\tabcolsep}{0.08cm} %donne la distance entre les collones du tableau%
\renewcommand{\arraystretch}{1} %donne la distance entre les lignes%
\begin{center}
{\begin{tabular}{ccccccccccccccc}
\hline
 & &  &  & & &  &  &  &  &  &  &   &     \\
 & $f(x)$  & method &  &    & $\lambda=0$ &  &  &$\lambda=0.25$ &  &  &  $\lambda=0.5$ & &   &  $\lambda=1$\\  

 &  &  & & &  &  & &  &   & &  &  &  &  \\
\hline
\hline
 & & T1&   &   & 0.040&   &  & 0.890  &  &  &1.000  & &   &1.000   \\
  & $x^2$& &   &    & &   &  &   &  &  &   & &   &   \\
  & &T2&   &   & 0.050&   &  & 1.000  &  &  &1.000  & &   &1.000  \\
\hline
 & &T1&   &   & 0.043&   &  & 1.000  &  &  &1.000  & &   &1.000   \\
  & $x^3$& &   &    & &   &  &   &  &  &   & &   &    \\
  & &T2&   &  & 0.045&   &  & 1.000  &  &  &1.000  & &   &1.000  \\
\hline
 & & T1&   &   & 0.035&   &  & 0.960  &  &  &1.000  & &   &1.000   \\
  & $x^2\cos(x)$& &   &    & &   &  &   &  &  &   & &   &   \\
  & & T2&   &   & 0.055&   &  & 1.000  &  &  &1.000  & &   &1.000   \\
\hline
 & & T1&   &   & 0.043&   &  & 0.036  &  &  &0.056  & &   &0.053 \\
  & $\sin(x)$& &   &    & &   &  &   &  &  &   & &   &    \\
  & & T2&   &   & 0.050&   &  & 0.065  &  &  &0.060  & &   &0.070   \\
\hline\hline
\end{tabular}}         
\end{center}}
\centering \caption{\label{table:tab4}Empirical sizes and powers over 300 replications  with significance level $\alpha=0.05$.}
\end{table}

\vskip 3mm

\noindent 5. PROOFS
\vskip 3mm
\noindent 5.1.  PRELIMINARY RESULTS

\vskip 3mm

\noindent We first give some technical lemmas which will be useful for proving the main theorems.

\begin{lem}\label{l1}
As $n\rightarrow +\infty$, we have
\begin{equation}\label{eqt3}
\begin{aligned}
 \sqrt{n}\left(\bigotimes_{\ell=1}^{d} \widehat{m}_\ell - \bigotimes_{\ell=1}^{d} m_\ell\right) \stackrel{\mathscr{D}}{\longrightarrow}\mathcal{W},
\end{aligned}
\end{equation}
where $\mathcal{W}$ is a random variable having the normal distribution in $\mathcal{H}$ with mean $0$ and covariance operator equal to that of 
\[
K_1 (X^{(1)},\cdot)\otimes\bigotimes_{\ell=2}^{d}m_\ell+ \sum_{\ell=2}^{d-1}\bigg(\bigotimes_{k=1}^{\ell-1} m_k \otimes K_\ell (X^{(\ell)},\cdot) \otimes \bigotimes_{k=\ell+1}^{d}m_k\bigg) + \bigotimes_{\ell=1}^{d-1} m_\ell \otimes K_d (X^{(d)},\cdot).
\]
\end{lem}
\noindent\textit{Proof}. Let us consider the random variables valued into $ \mathcal{H}_1 \times \mathcal{H}_2 \times \dots \times \mathcal{H}_d$ defined as   
\[
Z =\Big(K_1 (X^{(1)},\cdot),K_2 (X^{(2)},\cdot),\dots,K_d (X^{(d)},\cdot)\Big),
\] 
\[
Z_i =\Big(K_1 (X^{(1)}_i,\cdot),K_2 (X^{(2)}_i,\cdot),\dots,K_d (X^{(d)}_i,\cdot)\Big),\,\,i=1,\cdots,n,
\] 
with mean equal to 
$\mu =(m_1,m_2,\dots,m_d)$. For each $\ell\in\{1,\cdots,d\}$, we have $\sqrt{n}\left( \widehat{m}_\ell - m_\ell\right)=\pi_\ell(\widehat{H}_n)$, where  $\widehat{H}_n  = \sqrt{n} \Big ( \frac{1}{n}  \sum_{i=1}^{n} Z_i - \mu \Big )$ and $\pi_\ell$ is the canonical projection:
\begin{equation*}
\begin{aligned}
\pi_\ell : (f_1,f_2,\dots,f_d) \in \mathcal{H}_1 \times \mathcal{H}_2 \times \dots \times \mathcal{H}_d \longmapsto f_\ell \in \mathcal{H}_\ell.
\end{aligned}
\end{equation*}
Then, from the decomposition
\begin{equation*}
\begin{aligned}
\sqrt{n} \bigg(   \bigotimes_{\ell=1}^{d} \widehat{m}_\ell - \bigotimes_{\ell=1}^{d} m_\ell  \bigg) & = \bigg(\sqrt{n} (\widehat{m}_{1} - m_{1})\bigg)\otimes\bigotimes_{\ell=2}^{d}\widehat{m}_\ell\\
&+ \sum_{\ell=2}^{d-1}\bigg(\bigotimes_{k=1}^{\ell-1} m_k \otimes \left(\sqrt{n} (\widehat{m}_{\ell} - m_{\ell})\right) \otimes \bigotimes_{k=\ell+1}^{d}\widehat{m}_k\bigg) \\
&  \qquad + \bigotimes_{\ell=1}^{d-1} m_\ell \otimes \bigg(\sqrt{n} (\widehat{m}_{d} - m_{d})\bigg)
\end{aligned}
\end{equation*}
we get 
\begin{equation}\label{produits}
\sqrt{n} \bigg(   \bigotimes_{\ell=1}^{d} \widehat{m}_\ell - \bigotimes_{\ell=1}^{d} m_\ell  \bigg) = \widehat{\Phi}_n \big( \widehat{H}_n \big),
\end{equation} 
where $\widehat{\Phi}_n$ is the random operator from $ \mathcal{H}_1 \times \mathcal{H}_2 \times \dots \times \mathcal{H}_d$ to $ \mathcal{H}$ such that
\begin{equation*}\label{phin}
\widehat{\Phi}_n(T)=\pi_1(T)\otimes\bigotimes_{\ell=2}^{d}\widehat{m}_\ell+ \sum_{\ell=2}^{d-1}\bigg(\bigotimes_{k=1}^{\ell-1} m_k \otimes\pi_\ell(T) \otimes \bigotimes_{k=\ell+1}^{d}\widehat{m}_k\bigg) + \bigotimes_{\ell=1}^{d-1} m_\ell \otimes \pi_d(T).
\end{equation*}
Considering the operator $\Phi$ from 
 $ \mathcal{H}_1 \times \mathcal{H}_2 \times \dots \times \mathcal{H}_d$ to $ \mathcal{H}$ such that
\begin{equation}\label{phi}
\Phi(T)=\pi_1(T)\otimes\bigotimes_{\ell=2}^{d}m_\ell+ \sum_{\ell=2}^{d-1}\bigg(\bigotimes_{k=1}^{\ell-1} m_k \otimes\pi_\ell(T) \otimes \bigotimes_{k=\ell+1}^{d}m_k\bigg) + \bigotimes_{\ell=1}^{d-1} m_\ell \otimes \pi_d(T),
\end{equation}
we have
\begin{eqnarray}\label{ineg}
\bigg\Vert \Phi_n \big( \widehat{H}_n \big) - \Phi \big( \widehat{H}_n \big)\bigg\Vert_\mathcal{H} & \leq &\bigg\Vert\pi_1 \big( \widehat{H}_n \big) \otimes \bigg( \bigotimes_{\ell=2}^{d} \widehat{m}_\ell - \bigotimes_{\ell=2}^{d} m_\ell   \bigg)\bigg\Vert_\mathcal{H}\nonumber \\
& &+ \sum_{\ell=2}^{d-1}\bigg\Vert  \bigotimes_{p=1}^{\ell-1} m_\ell \otimes \pi_\ell \big( \widehat{H}_n \big) \otimes \bigg( \bigotimes_{p=\ell+1}^{d} \widehat{m}_p - \bigotimes_{p=\ell+1}^{d} m_p \bigg)  \bigg\Vert_\mathcal{H} \nonumber\\
 & \leq &\bigg\Vert\pi_1 \big( \widehat{H}_n \big) \bigg\Vert_\mathcal{H} \,\,  \bigg\Vert \bigotimes_{\ell=2}^{d} \widehat{m}_\ell - \bigotimes_{\ell=2}^{d} m_\ell  \bigg\Vert_\mathcal{H}\nonumber \\
&  +& \sum_{\ell=2}^{d-1}\bigg\Vert  \bigotimes_{p=1}^{\ell-1} m_\ell  \bigg\Vert_\mathcal{H}\, \bigg\Vert \pi_\ell \big( \widehat{H}_n \big) \bigg\Vert_\mathcal{H}\,  \bigg\Vert \bigotimes_{p=\ell+1}^{d} \widehat{m}_p - \bigotimes_{p=\ell+1}^{d} m_p   \bigg\Vert_\mathcal{H}.
\end{eqnarray}
The cental limit  theorem ensures that $\widehat{H}_n$ converges in distribution, as $n\rightarrow +\infty$, to a random variable $H$ having a centered normal distribution in $\mathcal{H}$ with covariance operator equal to that of $Z$. Since $\pi_\ell$ is continuous, we deduce that $\pi_\ell \big( \widehat{H}_n \big)$ converges in distribution, as $n\rightarrow +\infty$, to $\pi_\ell (H)$.  On the other hand, from the law of large numbers each $\widehat{m}_\ell$ converges almost surely, as $n\rightarrow +\infty$, to $m_\ell$. Then, for any $k\in\{1,\cdots,d-1\}$, we deduce from the inequality
\begin{equation*}
\begin{aligned}
\bigg\Vert   \bigotimes_{\ell=k}^{d} \widehat{m}_\ell - \bigotimes_{\ell=k}^{d} m_\ell  \bigg\Vert _\mathcal{H} & =\bigg\Vert ( \widehat{m}_{k} - m_{k})\otimes\bigotimes_{\ell=k+1}^{d}\widehat{m}_\ell\\
&+ \sum_{\ell=k+1}^{d-1}\bigg(\bigotimes_{j=1}^{\ell-1} m_j \otimes \left(\widehat{m}_{\ell} - m_{\ell}\right) \otimes \bigotimes_{j=\ell+1}^{d}\widehat{m}_j\bigg) \\
&  \qquad + \bigotimes_{\ell=1}^{d-1} m_\ell \otimes (\widehat{m}_{d} - m_{d}) \bigg\Vert _\mathcal{H} \\
& \leq \big\Vert \widehat{m}_{k} - m_{k}\big\Vert_{\mathcal{H}_k}\prod_{\ell=k+1}^{d} \big\Vert\widehat{m}_\ell \big\Vert_{\mathcal{H}_\ell}\\
&+ \sum_{\ell=k+1}^{d-1}\bigg(\prod_{j=1}^{\ell-1}\big\Vert m_j \big\Vert_{\mathcal{H}_j} \big\Vert \widehat{m}_{\ell} - m_{\ell}\big\Vert_{\mathcal{H}_\ell} \prod_{j=\ell+1}^{d}\big\Vert \widehat{m}_j\big\Vert_{\mathcal{H}_j}\bigg) \\
&  \qquad + \bigg(\prod_{\ell=1}^{d-1} \big\Vert m_\ell \big\Vert_{\mathcal{H}_\ell}\bigg)  \big\Vert\widehat{m}_{d} - m_{d} \big\Vert _{\mathcal{H}_d} 
\end{aligned}
\end{equation*}
that $\big\Vert   \bigotimes_{\ell=k}^{d} \widehat{m}_\ell - \bigotimes_{\ell=k}^{d} m_\ell  \big\Vert _\mathcal{H} $  converges almost surely, as $n\rightarrow +\infty$, to $0$. Therefore, from \eqref{ineg} we deduce that  $\Phi_n \big( \widehat{H}_n \big) - \Phi \big( \widehat{H}_n \big)$ converges in probability, as $n\rightarrow +\infty$, to $0$; consequently, $\sqrt{n} \bigg(   \bigotimes_{\ell=1}^{d} \widehat{m}_\ell - \bigotimes_{\ell=1}^{d} m_\ell  \bigg) = \Phi \big( \widehat{H}_n \big)+o_P(1)$  and, since $\Phi$ is contionuous,  Slustky's theorrem implies  the convergence in distribution  , as $n\rightarrow +\infty$, to $\mathcal{W}=\Phi (H)$. This random variable has a centered normal distribution in $\mathcal{H}$ with covariance operator equal to that of 
\begin{eqnarray*}
\Phi(Z)&=&
K_1 (X^{(1)},\cdot)\otimes\bigotimes_{\ell=2}^{d}m_\ell+ \sum_{\ell=2}^{d-1}\bigg(\bigotimes_{k=1}^{\ell-1} m_k \otimes K_\ell (X^{(\ell)},\cdot) \otimes \bigotimes_{k=\ell+1}^{d}m_k\bigg) \\
&&+ \bigotimes_{\ell=1}^{d-1} m_\ell \otimes K_d (X^{(d)},\cdot).
\end{eqnarray*}
\hfill $\Box$

\begin{lem}\label{l3}
Let $(Z_i)_{1 \leq i \leq n}$ be an i.i.d.  sample of a random variable  $Z$  valued into a Hilbert space $H$ and admitting a mean $m_Z$. If there exists a  real number $M>0$ such that $\Vert Z\Vert\leq  M$, then the statistic
\begin{equation}\label{eqt2}
\begin{aligned}
\widehat{s}_n  ^2= \frac{1}{n} \sum_{i=1}^{n} \big \langle Z_i , \overline{Z}_n \big \rangle^2 - \bigg(\frac{1}{n} \sum_{i=1}^{n} \big \langle Z_i , \overline{Z}_n \big \rangle\bigg)^2,
\end{aligned}
\end{equation}
where $\overline{Z}_n=\frac{1}{n} \sum_{i=1}^{n} Z_i$, is a consistent estimator of $Var \left(\big \langle Z , m_Z \big \rangle \right)$.
\end{lem}
\noindent\textit{Proof}. We have
$$  \frac{1}{n} \sum_{i=1}^{n} \big \langle Z_i , \overline{Z}_n \big \rangle^2 - \frac{1}{n} \sum_{i=1}^{n} \big \langle Z_i , m_Z \big \rangle^2  = \frac{1}{n} \sum_{i=1}^{n} \big \langle Z_i , \overline{Z}_n - m_Z\big \rangle^2 + \frac{2}{n} \sum_{i=1}^{n} \big \langle Z_i , m_Z \big \rangle \big \langle Z_i , \overline{Z}_n - m_Z \big \rangle ; $$
and using the Cauchy-Shwarz inequality, we get the inequalities
$$  \bigg|\frac{1}{n} \sum_{i=1}^{n} \big \langle Z_i , \overline{Z}_n - m_Z\big \rangle^2 \bigg\vert  \leq M^2 \big\Vert\overline{Z}_n - m_Z\big\Vert^2$$
and
$$ \bigg| \frac{2}{n} \sum_{i=1}^{n} \big \langle Z_i , m_Z \big \rangle \big \langle Z_i , \overline{Z}_n - m_Z \big \rangle \bigg|  \leq 2 M^2 \big\Vert m_Z\big\Vert\,\big\Vert \overline{Z}_n - m_Z\big\Vert$$
from which we conclude that 
$$\frac{1}{n} \sum_{i=1}^{n} \big \langle Z_i , \overline{Z}_n \big \rangle^2 - \frac{1}{n} \sum_{i=1}^{n} \big \langle Z_i , m_Z \big \rangle^2 = o_P(1)$$
since  $\big\Vert \overline{Z}_n - m_Z\big\Vert=o_P(1)$.
The law of large numbers and Slutsky’s theorem allow to conclude that the sequence $\frac{1}{n} \sum_{i=1}^{n} \big \langle Z_i , \overline{Z}_n \big \rangle^2$ converges in probability to $\mathbb{E}\left(  \big \langle Z , m_Z \big \rangle^2 \right)$ as $n \rightarrow +\infty$. Similarly, from
\[
\bigg|\frac{1}{n} \sum_{i=1}^{n} \big \langle Z_i , \overline{Z}_n \big \rangle - \frac{1}{n} \sum_{i=1}^{n} \big \langle Z_i , m_Z \big \rangle  \bigg|  = \bigg|\frac{1}{n} \sum_{i=1}^{n} \big \langle Z_i , \overline{Z}_n - m_Z\big \rangle \bigg|  \leq M \big\Vert \overline{Z}_n - m_Z\big\Vert,
\]
we get $\frac{1}{n} \sum_{i=1}^{n} \big \langle Z_i , \overline{Z}_n \big \rangle - \frac{1}{n} \sum_{i=1}^{n} \big \langle Z_i , m_Z \big \rangle = o_P(1)$, therefore,  and $\frac{1}{n} \sum_{i=1}^{n} \big \langle Z_i , \overline{Z}_n \big \rangle$ converges in probability to $\mathbb{E}\left(  \big \langle Z , m_Z \big \rangle \right)$ as $n\rightarrow +\infty$ just like  $ \frac{1}{n} \sum_{i=1}^{n} \big \langle Z_i , m_Z \big \rangle$. So, $\widehat{s}_n ^2$ is a consistent estimator of $Var \left( \big \langle Z , m_Z \big \rangle\right)$.
\hfill $\Box$

\vskip 3mm
\noindent 5.2.  PROOF OF THEOREM \ref{thm}

\vskip 3mm

\noindent Using  \eqref{produit} and  the reproducing properties of $K$, it is easy to see that
\begin{eqnarray}
\widehat{\textit{\textbf{D}}}_{\gamma,n}&=&\bigg\Vert\frac{1}{n}\sum_{i=1}^n  K(X_i,\cdot) \bigg\Vert_{\mathcal{H}} ^2 +\bigg\Vert \bigotimes_{\ell=1}^{d} \widehat{m}_\ell \bigg\Vert_{\mathcal{H}} ^2- \frac{2}{n} \sum_{i=1}^{n}w_{i,n}(\gamma)\Big\langle  K(X_i,\cdot),   \bigotimes_{\ell=1}^{d} \widehat{m}_\ell  \Big\rangle_{\mathcal{H}}\nonumber \\
&=&\big\Vert\widehat{m}-m\big\Vert_{\mathcal{H}} ^2 -\big\Vert m\big\Vert_{\mathcal{H}} ^2 +2\langle \widehat{m},m\rangle_\mathcal{H}+\bigg\Vert \bigotimes_{\ell=1}^{d} \widehat{m}_\ell - \bigotimes_{\ell=1}^{d}  m_\ell \bigg\Vert_{\mathcal{H}} ^2-\big\Vert  \bigotimes_{\ell=1}^{d} m_\ell \big\Vert _\mathcal{H} ^2 \nonumber\\
& &+2\Big\langle  \bigotimes_{\ell=1}^{d}  \widehat{m}_\ell ,  \bigotimes_{\ell=1}^{d}  m_\ell\Big\rangle_\mathcal{H} - \frac{2}{n} \sum_{i=1}^{n}w_{i,n}(\gamma)\Big\langle  K(X_i,\cdot),   \bigotimes_{\ell=1}^{d} \widehat{m}_\ell  \Big\rangle_{\mathcal{H}}\nonumber ,
\end{eqnarray}
where $\widehat{m}=\frac{1}{n}\sum_{i=1}^n  K(X_i,\cdot) $. Since $d\textrm{HSIC}(X) 
 =\big \Vert   m \big\Vert _\mathcal{H} ^2+\big\Vert  \bigotimes_{\ell=1}^{d} m_\ell \big\Vert _\mathcal{H} ^2+2\Big\langle  m,  \bigotimes_{\ell=1}^{d}  m_\ell\Big\rangle_\mathcal{H}$ and 
\begin{eqnarray*}
\frac{1}{n} \sum_{i=1}^{n} w_{i,n}(\gamma)\Big \langle   K(X_i,\cdot) ,  \bigotimes_{\ell=1}^{d} \widehat{m}_\ell  \Big \rangle_\mathcal{H}  
& = &	\Big \langle  \frac{1}{n} \sum_{i=1}^{n} \left(w_{i,n}(\gamma) - 1 \right) K(X_i,\cdot) ,  \bigotimes_{\ell=1}^{d} \widehat{m}_\ell - \bigotimes_{\ell=1}^{d} m_\ell  \Big \rangle_\mathcal{H}  \\
& &	+ \Big \langle  \frac{1}{n} \sum_{i=1}^{n} w_{i,n}(\gamma) K(X_i,\cdot) ,  \bigotimes_{\ell=1}^{d} m_\ell  \Big \rangle_\mathcal{H}\\
& &+ \Big \langle   \widehat{m} - m,  \bigotimes_{\ell=1}^{d} \widehat{m}_\ell - \bigotimes_{\ell=1}^{d} m_\ell  \Big \rangle_\mathcal{H}  \\
&  &+ \Big \langle  m,  \bigotimes_{\ell=1}^{d} \widehat{m}_\ell \Big \rangle_\mathcal{H}  -  \frac{1}{n} \sum_{i=1}^{n} w_{i,n}(\gamma) \Big \langle m, \bigotimes_{\ell=1}^{d} m_\ell  \Big \rangle_\mathcal{H} \\
& &+ \frac{1}{n} \sum_{i=1}^{n} \left(w_{i,n}(\gamma) - 1 \right)\Big \langle  m,\bigotimes_{\ell=1}^{d} m_\ell  \Big \rangle_\mathcal{H} , 
\end{eqnarray*}
it follows $\sqrt{n}\bigg(\widehat{\textit{\textbf{D}}}_{\gamma,n}-d\textrm{HSIC}(X)\bigg)=A_n+B_n+C_n$, where
\begin{eqnarray*}
A_n 
& = &\frac{1}{\sqrt{n}}   \Vert    \sqrt{n}\left( \widehat{m}-m \right)\Vert _\mathcal{H}^2 
+\frac{1}{\sqrt{n}} \bigg\Vert  \sqrt{n} \bigg(\bigotimes_{\ell=1}^{d} \widehat{m}_\ell - \bigotimes_{\ell=1}^{d} m_\ell \bigg)\bigg\Vert _\mathcal{H}^2 \\
&&-2  \Big \langle   \widehat{m}-m ,  \sqrt{n} \bigg(\bigotimes_{\ell=1}^{d} \widehat{m}_\ell - \bigotimes_{\ell=1}^{d} m_\ell \bigg)  \Big \rangle_\mathcal{H},
\end{eqnarray*}
\begin{eqnarray*}
B_n  & =& -2  \Big \langle  \frac{1}{n} \sum_{i=1}^{n} \left(w_{i,n}(\gamma) - 1 \right) K(X_i,\cdot) ,  \sqrt{n} \bigg(\bigotimes_{\ell=1}^{d} \widehat{m}_\ell - \bigotimes_{\ell=1}^{d} m_\ell \bigg)  \Big \rangle_\mathcal{H} \\
& &- 2 \sqrt{n} \bigg( \frac{1}{n} \sum_{i=1}^{n} \left(w_{i,n}(\gamma) - 1 \right)\bigg) \Big \langle  m, \bigotimes_{\ell=1}^{d} m_\ell   \Big \rangle_\mathcal{H}
\end{eqnarray*}
and
\begin{eqnarray}\label{cn}
C_n & =& \sqrt{n} \bigg(	 - \Vert   m \Vert _\mathcal{H} ^2 +2   \Big \langle\widehat{m} ,m \Big \rangle_{\mathcal{H}}  - \Vert    \bigotimes_{\ell=1}^{d} m_\ell  \Vert _\mathcal{H} ^2 + 2 \Big \langle  \bigotimes_{\ell=1}^{d}  \widehat{m}_\ell  ,  \bigotimes_{\ell=1}^{d} m_\ell  \Big \rangle_{\mathcal{H}} \nonumber\\
&  &	 	- 2\Big \langle  \frac{1}{n} \sum_{i=1}^{n} w_{i,n}(\gamma) K(X_i,\cdot) ,   \bigotimes_{\ell=1}^{d} m_\ell   \Big \rangle_\mathcal{H}	- 2 	\Big \langle  m,   \bigotimes_{\ell=1}^{d}  \widehat{m}_\ell  \Big \rangle_\mathcal{H} \nonumber\\
&  &+  \frac{2}{n} \sum_{i=1}^{n} w_{i,n}(\gamma) \Big \langle m,  \bigotimes_{\ell=1}^{d} m_\ell   \Big \rangle_\mathcal{H}		 - \Vert    m\Vert _\mathcal{H} ^2 - \Vert    \bigotimes_{\ell=1}^{d} m_\ell  \Vert _\mathcal{H} ^2 + 2\Big \langle  m , \bigotimes_{\ell=1}^{d} m_\ell  \Big \rangle_\mathcal{H}\bigg)\nonumber\\
& =&\sqrt{n} \bigg(		2 \Big \langle	\widehat{m} - m , m  \Big \rangle_\mathcal{H} - \frac{2}{n} \sum_{i=1}^{n} w_{i,n}(\gamma) \Big \langle	 K(X_i,\cdot) - m, \bigotimes_{\ell=1}^{d} m_\ell  \Big \rangle_\mathcal{H} \nonumber\\
& & -2 \Big \langle	  \bigotimes_{\ell=1}^{d}  \widehat{m}_\ell  -  \bigotimes_{\ell=1}^{d} m_\ell   ,m -  \bigotimes_{\ell=1}^{d} m_\ell  \Big \rangle_\mathcal{H} \bigg).
\end{eqnarray}
First, using Cauchy-Schwartz inequality, we get
\begin{eqnarray*}
\vert A_n\vert&\leq&\frac{1}{\sqrt{n}}   \Vert    \sqrt{n}\left( \widehat{m}-m \right)\Vert _\mathcal{H}^2 
+\frac{1}{\sqrt{n}} \bigg\Vert  \sqrt{n} \bigg(\bigotimes_{\ell=1}^{d} \widehat{m}_\ell - \bigotimes_{\ell=1}^{d} m_\ell \bigg)\bigg\Vert _\mathcal{H}^2 \\
&&+2   \Vert     \widehat{m}-m  \Vert _\mathcal{H} \bigg\Vert  \sqrt{n} \bigg(\bigotimes_{\ell=1}^{d} \widehat{m}_\ell - \bigotimes_{\ell=1}^{d} m_\ell \bigg)\bigg\Vert _\mathcal{H}.
\end{eqnarray*}
The law of large numbers  ensures that $ \Vert     \widehat{m}-m  \Vert _\mathcal{H}\rightarrow 0$ as $n\rightarrow +\infty$. Then, from Lemma \ref{l1} we deduce that $\frac{1}{\sqrt{n}} \bigg\Vert  \sqrt{n} \bigg(\bigotimes_{\ell=1}^{d} \widehat{m}_\ell - \bigotimes_{\ell=1}^{d} m_\ell \bigg)\bigg\Vert _\mathcal{H}^2 $ and $\Vert     \widehat{m}-m  \Vert _\mathcal{H} \bigg\Vert  \sqrt{n} \bigg(\bigotimes_{\ell=1}^{d} \widehat{m}_\ell - \bigotimes_{\ell=1}^{d} m_\ell \bigg)\bigg\Vert _\mathcal{H}$ converge in probability to $0$ as $n\rightarrow +\infty$. The central limit theorem gives the convergence in distrubution of  $\widehat{m}-m$, as  $n\rightarrow +\infty$, to a random variable having a normal distribution in $\mathcal{H}$, hence  $\frac{1}{\sqrt{n}}   \Vert    \sqrt{n}\left( \widehat{m}-m \right)\Vert _\mathcal{H}^2$ converges in probability to $0$ as $n\rightarrow +\infty$. Then, from the preceding inequality we get $A_n=o_P(1)$. Secondly, another use of Cauchy-Schwartz inequality gives for $n$ large enough
\begin{eqnarray}\label{bn}
\vert B_n\vert &\leq & \bigg \Vert  \frac{1}{n} \sum_{i=1}^{n} \left(w_{i,n}(\gamma) - 1 \right) K(X_i,\cdot)\bigg\Vert _\mathcal{H} \,\, \bigg \Vert \sqrt{n} \bigg(  \bigotimes_{\ell=1}^{d}  \widehat{m}_\ell  -  \bigotimes_{\ell=1}^{d} m_\ell  \bigg)\bigg \Vert _\mathcal{H}\nonumber \\
&  &\qquad \qquad + 2\sqrt{n} \Big | \frac{1}{n} \sum_{i=1}^{n} \left(w_{i,n}(\gamma) - 1 \right)\Big |\,\, \Vert   m \Vert _\mathcal{H} \,\,\bigg\Vert     \bigotimes_{\ell=1}^{d} m_\ell   \bigg\Vert _\mathcal{H}\nonumber\\
&\leq & \bigg \Vert  \frac{1}{n} \sum_{i=1}^{n} \left(w_{i,n}(\gamma) - 1 \right) K(X_i,\cdot)\bigg\Vert _\mathcal{H} \,\, \bigg \Vert \sqrt{n} \bigg(  \bigotimes_{\ell=1}^{d}  \widehat{m}_\ell  -  \bigotimes_{\ell=1}^{d} m_\ell  \bigg)\bigg \Vert _\mathcal{H}\nonumber \\
&  &\qquad \qquad +\frac{ 2\tau}{\sqrt{n}}  \Vert   m \Vert _\mathcal{H} \,\,\bigg\Vert     \bigotimes_{\ell=1}^{d} m_\ell   \bigg\Vert _\mathcal{H}
\end{eqnarray}
since, from $(\mathscr{A}_2)$, we have $n \Big | \frac{1}{n} \sum_{i=1}^{n} \left(w_{i,n}(\gamma) - 1 \right)\Big |\leq \tau$. Lemma 1 in Manfoumbi Djonguet et al. (2022) ensures that $ \Vert  \frac{1}{n} \sum_{i=1}^{n} \left(w_{i,n}(\gamma) - 1 \right) K(X_i,\cdot) \Vert _\mathcal{H}=o_P(1)$, and since the sequence  $ \Vert \sqrt{n} (  \bigotimes_{\ell=1}^{d}  \widehat{m}_\ell  -  \bigotimes_{\ell=1}^{d} m_\ell  )  \Vert _\mathcal{H}$ converges in distribution, as $n\rightarrow +\infty$, to an appropriate random variable (see Lemma 1), we deduce from \eqref{bn}  that $B_n=o_P(1)$. Thirdly, 
reporting \eqref{produits}  in \eqref{cn} gives
\begin{eqnarray}\label{cn2}
C_n & = &\sqrt n \Big \{	\frac{2}{n} \sum_{i=1}^{n} \Big \langle	 K(X_i,\cdot) - m   , m   \Big \rangle_ \mathcal{H} - \frac{2}{n} \sum_{i=1}^{n} w_{i,n}(\gamma) \Big \langle	 K(X_i,\cdot) - m   , \bigotimes_{\ell=1}^{d} m_\ell \Big \rangle_ \mathcal{H} \Big \} \nonumber\\
& & \qquad  + 2 \Big \langle	\Phi_n \big( \widehat{H}_n \big) - \Phi \big( \widehat{H}_n \big)   ,\bigotimes_{\ell=1}^{d} m_\ell - m  \Big \rangle_ \mathcal{H} + 2 \Big \langle \Phi \big( \widehat{H}_n \big)   , \bigotimes_{\ell=1}^{d} m_\ell - m  \Big \rangle_ \mathcal{H}.
\end{eqnarray}
From
\[
\bigg\vert \Big \langle	\Phi_n \big( \widehat{H}_n \big) - \Phi \big( \widehat{H}_n \big)   ,\bigotimes_{\ell=1}^{d} m_\ell - m  \Big \rangle_ \mathcal{H}\bigg\vert\leq \bigg\Vert 	\Phi_n \big( \widehat{H}_n \big) - \Phi \big( \widehat{H}_n \big) \bigg\Vert_ \mathcal{H} \bigg\Vert\bigotimes_{\ell=1}^{d} m_\ell - m  \bigg\Vert_ \mathcal{H}
\]
and the convergence in probability of $	\Phi_n \big( \widehat{H}_n \big) - \Phi \big( \widehat{H}_n \big) $  to $0$ as $n\rightarrow +\infty$  (see Lemma 1), we deduce that 
\begin{equation}\label{op}
\Big \langle	\Phi_n \big( \widehat{H}_n \big) - \Phi \big( \widehat{H}_n \big)   ,\bigotimes_{\ell=1}^{d} m_\ell - m  \Big \rangle_ \mathcal{H} = o_P(1).
\end{equation}
We will now decompose the term $\Big \langle \Phi \big( \widehat{H}_n \big)   , \bigotimes_{\ell=1}^{d} m_\ell - m  \Big \rangle_ \mathcal{H}$ using (\ref{phi}).  For ease of notation we  set   $\nu_{k}^q = \bigotimes_{\ell=k}^{q} m_\ell$, and we have
\begin{equation*}
\begin{aligned}
\pi_1 \big( \widehat{H}_n \big)\otimes \bigotimes_{\ell=2}^{d} m_\ell & =  \sqrt{n} (\hat{m}_1 - m_1) \otimes \nu_{2}^{q} \\
& = \frac{1}{\sqrt{n}}\sum_{i=1}^{n}\left(K_1(X_i ^{(1)},\cdot) - m_1\right) \otimes \nu_2^{d}\\
& = \frac{1}{\sqrt{n}}\sum_{i=1}^{n}\bigg(  K_1(X_i ^{(1)},\cdot)\otimes \nu_{2}^{d} -  \bigotimes_{\ell=1}^{d} m_\ell  \bigg),\\
\end{aligned}
\end{equation*}
\begin{equation*}
\begin{aligned}
\bigotimes_{k=1}^{\ell-1} m_k \otimes \pi_\ell \big( \widehat{H}_n \big) \otimes \bigotimes_{k=\ell+1}^{d} m_k & = \nu_{1}^{\ell-1} \otimes \sqrt{n} (\widehat{m}_\ell - m_\ell) \otimes \nu_{\ell+1}^{d} \\
& = \frac{1}{\sqrt{n}}\sum_{i=1}^{n}\nu_{1}^{\ell-1} \otimes \left(K_\ell(X_i ^{(\ell)},\cdot) - m_l\right) \otimes \nu_{\ell+1}^{d}\\
& = \frac{1}{\sqrt{n}}\sum_{i=1}^{n}\bigg(\nu_{1}^{\ell-1} \otimes K_\ell(X_i ^{(\ell)},\cdot)\otimes \nu_{\ell+1}^{d} -  \bigotimes_{\ell=1}^{d} m_\ell  \bigg),\\
\end{aligned}
\end{equation*}
and
\begin{equation*}
\begin{aligned}
\bigotimes_{\ell=1}^{d-1} m_\ell \otimes \pi_d \big( \widehat{H}_n \big) & = \nu_{1}^{d-1} \otimes \left(\sqrt{n} (\widehat{m}_d - m_d)\right)  \\
& = \frac{1}{\sqrt{n}}\sum_{i=1}^{n}\nu_{1}^{d-1} \otimes \left(K_d(X_i ^{(d)},\cdot) - m_d\right) \\
& = \frac{1}{\sqrt{n}}\sum_{i=1}^{n}\bigg(\nu_{1}^{d-1} \otimes K_d(X_i ^{(d)},\cdot) -  \bigotimes_{\ell=1}^{d} m_\ell  \bigg).\\
\end{aligned}
\end{equation*}
Hence
\[
\Phi \big( \widehat{H}_n \big) 
 = \frac{1}{\sqrt{n}}\sum_{i=1}^{n} \Big \{ K_1(X_i ^{(1)},\cdot)\otimes \nu_{2}^{d} + \sum_{\ell=2}^{d-1}   \left(\nu_{1}^{\ell-1} \otimes K_\ell(X_i ^{(\ell)},\cdot)\otimes \nu_{\ell+1}^{d}\right) + \nu_{1}^{d-1}\otimes K_d(X_i ^{(d)},\cdot) -  d\bigotimes_{\ell=1}^{d}m_\ell\Big \},
\]
and from \eqref{cn2} and \eqref{op}  it follows that $C_n=D_n+o_P(1)$, where
\begin{eqnarray*}
& &D_n   =  \frac{2}{\sqrt{n}}  \sum_{i=1}^{n}\bigg\{ \Big \langle	 K(X_i,\cdot) - m   , m   \Big \rangle_ \mathcal{H} -  w_{i,n}(\gamma) \Big \langle	 K(X_i,\cdot) - m   , \bigotimes_{\ell=1}^{d} m_\ell \Big \rangle_ \mathcal{H} \nonumber\\
&&+\Big\langle  K_1(X_i ^{(1)},\cdot)\otimes \nu_{2}^{d} + \sum_{\ell=2}^{d-1}   \left(\nu_{1}^{\ell-1} \otimes K_\ell(X_i ^{(\ell)},\cdot)\otimes \nu_{\ell+1}^{d}\right) \\
&&+ \nu_{1}^{d-1}\otimes K_d(X_i ^{(d)},\cdot) -  d\bigotimes_{\ell=1}^{d}m_\ell, \bigotimes_{\ell=1}^{d} m_\ell - m  \Big \rangle_ \mathcal{H}\bigg\}  +o_P(1)\\
&&= \frac{2}{\sqrt{n}}\sum_{i=1}^{n}\bigg( \mathcal{U} (X_i) - w_{i,n}(\gamma)\mathcal{V}  (X_i) \bigg)  + o_P (1) .
\end{eqnarray*}
Finally,   $\sqrt{n}\bigg(\widehat{\textit{\textbf{D}}}_{\gamma,n}-d\textrm{HSIC}(X)\bigg)=E_n+o_P(1)$, where $E_n=\frac{2}{\sqrt{n}}\sum_{i=1}^{n}\bigg( \mathcal{U} (X_i) - w_{i,n}(\gamma)\mathcal{V}  (X_i) \bigg)$
and, consequently, $\sqrt{n}\bigg(\widehat{\textit{\textbf{D}}}_{\gamma,n}-d\textrm{HSIC}(X)\bigg)$ has the same limiting distribution than $E_n$; it remains to derive this latter.  Let us set  
\[
s_{n,\gamma}^2=\sum_{i=1}^{n}Var\Big(\mathcal{U}(X_i)-w_{i,n}(\gamma)\mathcal{V}(X_i)\Big);
\] 
 by similar arguments as in the proof of Theorem 1 in Makigusa and Naito (2020) we obtain that, for any $\varepsilon>0$,
\[
s_{n,\gamma}^{-2}\sum_{i=1}^{n}\int_{\{x:| \mathcal{U}(x)-w_{i,n}(\gamma)\mathcal{V}(x) |>\varepsilon s_{n,\gamma}\}}^{} \bigg(\mathcal{U}(x)-w_{i,n}(\gamma)\mathcal{V}(x)\bigg)^2\,d\mathbb{P}_{X}(x) 
\]
converges to $0$ as $n\rightarrow +\infty$. Then, by Section 1.9.3 in Serfling (1980),   
$\sqrt{n}s_{n,\gamma}^{-1}\frac{E_n}{2}	\stackrel{\mathscr{D}}{\rightarrow} \mathcal{N}\left(0,1\right)$. However,
\begin{eqnarray*}
\left(\frac{s_{n,\gamma}}{\sqrt{n}}\right)^2
&=&Var\left(\mathcal{U}(X_1)\right)+\left(\frac{1}{n}\sum_{i=1}^{n}w^2_{i,n}(\gamma)\right)Var\left(\mathcal{V}(X_1)\right)\nonumber\\
&&-2\left(\frac{1}{n}\sum_{i=1}^{n}w_{i,n}(\gamma)\right)Cov\left(\mathcal{U}(X_1),\mathcal{V}(X_1)\right),
\end{eqnarray*}
from  $ (\mathscr{A}_2) $ and $ (\mathscr{A}_4) $, 
\[
\lim_{n\rightarrow +\infty}\left(\frac{1}{n}\sum_{i=1}^{n}w^2_{i,n}(\gamma)\right)=w^2(\gamma)\,\,\,\textrm{ and }\,\,\,\lim_{n\rightarrow +\infty}\left(\frac{1}{n}\sum_{i=1}^{n}w_{i,n}(\gamma)\right)=1.
\]
Thus
\[
\lim\limits_{n\rightarrow +\infty}\left(n^{-1}s_{n,\gamma}^2\right)=Var\left(\mathcal{U}(X_1)\right)+ w^2(\gamma)Var\left(\mathcal{V}(X_1)\right) -2Cov\left(\mathcal{U}(X_1),\mathcal{V}(X_1)\right)
\]
and, therefore, $E_n	\stackrel{\mathscr{D}}{\rightarrow} \mathcal{N}\left(0,\sigma_\gamma^2\right)$.

\vskip 3mm
\noindent 5.3.  PROOF OF PROPOSITION  \ref{prp}

\vskip 3mm

\noindent Putting  $Z = K (X,\cdot)$,  $Z_i = K (X_i,\cdot)$  and  $\overline{Z}_n=\frac{1}{n}\sum_{i=1}^nZ_i$, we have
\[
\Vert Z\Vert_\mathcal{H}=\sqrt{K(X,X)}=\sqrt{\prod_{\ell=1}^{d}K_\ell(X^{(\ell)},X^{(\ell)})}\leq \prod_{\ell=1}^{d}\Vert K_\ell\Vert_\infty^{1/2} ,
\]
\[
\big \langle Z_i , \overline{Z}_n \big \rangle_\mathcal{H}  = \frac{1}{n}\sum_{j=1}^{n} \big \langle Z_i , {Z}_j \big \rangle_\mathcal{H}  = \frac{1}{n}\sum_{j=1}^{n}  K (X_i,X_j)= \frac{1}{n}\sum_{j=1}^{n}  \Bigg[\prod_{\ell=1}^{d}K_\ell(X^{(\ell)}_i,X^{(\ell)}_j)\Bigg]
\]
and
\begin{eqnarray*}
 \frac{1}{n} \sum_{i=1}^{n} \big \langle Z_i , \overline{Z}_n \big \rangle^2 - \bigg(\frac{1}{n} \sum_{i=1}^{n} \big \langle Z_i , \overline{Z}_n \big \rangle\bigg)^2  
& = &\frac{1}{n}\sum_{i=1}^{n}  \Bigg ( \frac{1}{n}\sum_{j=1}^{n}  \Bigg[\prod_{\ell=1}^{d}K_\ell(X^{(\ell)}_i,X^{(\ell)}_j)\Bigg]\Bigg)^2 \\
&&\qquad-\bigg(\frac{1}{n^2} \sum_{i=1}^{n} \sum_{j=1}^{n}  \Bigg[\prod_{\ell=1}^{d}K_\ell(X^{(\ell)}_i,X^{(\ell)}_j)\Bigg]\bigg)^2\\
&=&\widehat{\alpha}_n^2.
\end{eqnarray*}
By applying Lemma \ref{l3} we obtain the convergence in probability of $ \widehat{\alpha}_n^2$ to $Var\left(\left\langle K(X,\cdot),m\right\rangle_\mathcal{H}\right)=Var\left(\mathcal{V} (X_1)\right)$  as $n\rightarrow +\infty$  and, therefore, that of $\widehat{\sigma}_\gamma^2$ to $\sigma_\gamma^2$.

\vskip 3mm
\noindent BIBLIOGRAPHY
\vskip 3mm

\noindent  Aghoukeng Jiofack, J.G., and G.M. Nkiet.   2010.  Testing for lack of dependence between functional variables.\textit{Statistics and  Probability  Letters} 80, 1210-1217.

\noindent Ahmad, I.A. 1993.  Modification of some goodness-of-fit statistics to yield
asymptotic normal null distribution. \textit{Biometrika}  80, 466–472.

\noindent Balogoun, A.K.S.,  G.M. Nkiet, and C.  Ogouyandjou. 2021.  Asymptotic normality of a generalized maximum mean discrepancy estimator. \textit{Statistics and  Probability  Letters} 169, 108961.

\noindent Balogoun, A.K.S.,  G.M. Nkiet, and C.  Ogouyandjou. 2022.  A $k$-sample test for functional data based on generalized maximum mean discrepancy. \textit{Lithuanian     Mathematical  Journal}, DOI:10.1007/10986-022-09572-x.

\noindent Berlinet,  A., and C. Thomas-Agnan. 2004.  \textit{Reproducing kernel Hilbert spaces in probability and statistics}. Kluwer.

\noindent G\'orecki, T.,  M. Krzysko,  and W. Wołynski. 2020. Independence test and canonical correlation analysis based on the alignment between kernel matrices for multivariate functional data. \textit{Artificial  Intelligence  Review} 53, 475--499.

\noindent Gretton, A.,  K.M. Borgwardt,  M.J.   Rasch, B.  Sch$\ddot{\textrm{o}}$lkopf, and A.J.  Smola.  2012. A kernel two-sample test. \textit{Journal   Machine  Learning   Research} 13, 723--776.

\noindent Gretton, A.,  O. Bousquet,  A.  Smola,   and  B. Sch$\ddot{\textrm{o}}$lkopf. 2005. Measuring statistical dependence with
Hilbert-Schmidt norms. In \textit{Algorithmic Learning Theory}  (eds S. Jain, H.U. Simon and E. Tomita), 63--77. Berlin: Springer.

\noindent Gretton, A.,  K. Fukumizu,  C.H.   Teo, L. Song,  B.  Sch$\ddot{\textrm{o}}$lkopf, and A.J.  Smola.  2007. A kernel statistical test of independence. In \textit{Advances in Neural Information Processing Systems} 20 (eds J.C. Platt, D. Koller, Y. Singer and S.T. Roweis), 585--592. New York: Curran Associates.

\noindent Kokoszka, P., I. Maslova,  J.  Sojka, and L.  Zhu.  2008. Testing for lack of dependence in the functional linear model. \textit{Canadian  Journal of   Statistics}  36, 1-16.

\noindent Lai, T., Z. Zhang,  Y. Wang, and L.  Kong. 2021. Testing independence of functional variables by angle covariance. \textit{Journal of  Multivariate  Analysis} 182, 104711.

\noindent Laumann, F., J. Von K$\ddot{\textrm{u}}$gelgen,   and M.  Brahona. 2021. Kernel two-sample and independence tests for non-stationary random processes. \textit{Engineering  Proceedings} 5, 31.

\noindent Makigusa, N., and  K.  Naito. 2020. Asymptotic normality of a consistent estimator of maximum mean discrepancy in Hilbert space. \textit{Statistics and  Probability  Letters} 156, 108596.

\noindent Manfoumbi Djonguet, T.K., A.  Mbina Mbina,  and  G.M. Nkiet. 2022. Testing independence of functional variables by an Hilbert-Schmidt independence criterion estimator. ArXiv:2206.11607.

\noindent Pfister, N., P. B$\ddot{\textrm{u}}$hlmann, B.  Sch$\ddot{\textrm{o}}$lkopf,   and  J. Peters. 2018. Kernel-based tests for joint independence. \textit{Journal of the  Royal  Statistical  Society  Series  B Statistical  Methodology} 80,   5–31.

\noindent Serfling, R.J. 1980. \textit{Approximation Theorems of Mathematical Statistics}. Wiley, New-York.

\end{document}